\documentclass[12pt]{article}
\usepackage[backref=page,colorlinks,bookmarksopen=true]{hyperref}
\usepackage{pdfsync}
\usepackage{amssymb}
\usepackage{latexsym}
\usepackage{amsmath,mathabx}
\usepackage[all]{xy}
\newtheorem{theorem}{Theorem}[section]
\newtheorem{lemma}[theorem]{Lemma}
\newtheorem{proposition}[theorem]{Proposition}
\newtheorem{corollary}[theorem]{Corollary}
\newtheorem{definition}[theorem]{Definition}
\newtheorem{example}[theorem]{Example}

\newtheorem{remark}[theorem]{Remark}

\newtheorem{notations}[theorem]{Notations}

\begin{document}
\title{Classifying (Weak) Coideal Subalgebras of Weak Hopf $C^*$-Algebras}
\author{Leonid Vainerman Jean-Michel Vallin \\
\\
\it {Dedicated to the Memory of Etienne Blanchard}}
\date{}
\maketitle

\begin{abstract} We develop a general approach to the problem of classification of weak coideal
$C^*$-subalgebras of weak Hopf $C^*$-algebras. As an example, we consider weak Hopf $C^*$-algebras
and their weak coideal $C^*$-subalgebras associated with Tambara Yamagami categories.
\end{abstract}
\tableofcontents

\footnote {AMS Subject Classification [2010]{: Primary 16T05, Secondary 18D10, Tertiary
46L05.}}

\footnote{ Keywords {: Coactions and corepresentations
of quantum groupoids, $C^*$-categories, reconstruction theorem.}}

\newpage
\newenvironment{dm}{\hspace*{0,15in} {{\bf Proof.}}}{$\square$}
\begin{section}{Introduction}

It is known that any finite tensor category equipped with a fiber functor to the category of finite dimensional vector spaces is equivalent to the representation category of some Hopf algebra - see, for example, \cite{EGNO}, Theorem 5.3.12. But many tensor categories do not admit a fiber functor, so they cannot be presented as representation categories of Hopf algebras. On the other hand, T. Hayashi \cite{Ha} showed that any fusion category always admits a tensor functor to the category of bimodules over some semisimple (even commutative) algebra. Using this, it was proved in \cite{Ha}, \cite{Os}, \cite{Sz} that any fusion category is equivalent to the representation category of some algebraic structure generalizing Hopf algebras called a {\it weak Hopf algebra} \cite{BNSz} or a {\it finite quantum groupoid} \cite{NV}. The main difference between weak and usual Hopf algebra is that in the former the coproduct $\Delta$ is not necessarily unital.

Apart from tensor categories, weak Hopf algebras have interesting applications to the subfactor theory. In particular, for any finite index and finite depth $II_1$-subfactor $N\subset M$, there exists a weak Hopf $C^*$-algebra $\mathfrak G$ such that the corresponding Jones tower can be expressed in terms of crossed products of $N$ and $M$ with $\mathfrak G$ and its dual. Moreover, there is {\it a Galois correspondence} between intermediate subfactors in this Jones tower and coideal $C^*$-subalgebras of $\mathfrak G$ - see \cite{NV2}. This motivates the study of coideal $C^*$-subalgebras
of weak Hopf $C^*$-algebras which is the subject of the present paper. The abbreviation WHA will always mean a weak Hopf $C^*$-algebra.

A coideal $C^*$-subalgebra is a special case of the notion of a $\mathfrak G$-$C^*$-algebra, which is, by definition, a unital $C^*$ algebra $A$ equipped with a {\it coaction} $\mathfrak a$ of a WHA $\mathfrak G=(B,\Delta, S,\varepsilon)$. More exactly, we will use the following

\begin{definition} \label{weak} A weak right coideal $C^*$-subalgebra of $B$ is a right $\mathfrak G$-$C^*$-algebra $(A,\mathfrak a)$ with a $C^*$-algebra inclusion $i: A\mapsto B$ (not necessarily unital) satisfying $\Delta=
(i\otimes id_B)\mathfrak a$. One can think of $A$ as of a $C^*$-subalgebra of $B$ such that $\mathfrak a=\Delta$.
If $i$ is unital, we call $A$ a coideal $C^*$-subalgebra of $B$.
\end{definition}

For the sake of brevity, we will call a (weak) coideal $C^*$-subalgebra a (weak) coideal of $B$. Note that if $\mathfrak G$
is a usual Hopf $C^*$-algebra, then one can prove that necessarily $1_A=1_B$, so weak and usual coideals coincide.

It was shown in \cite{VV1} that any $\mathfrak G$-$C^*$-algebra $(A,\mathfrak a)$ corresponds to a pair $(\mathcal M,M)$, where $\mathcal M$ is a module category with a generator $M$ over the category of unitary corepresentations of $\mathfrak G$.

In Preliminaries we recall definitions and facts needed for the exact formulation of this result  expressed in Theorem \ref{main}. Note that similar categorical duality for compact quantum group coactions was obtained earlier in \cite{CY}, \cite{Nes1}.  

Section 3  is devoted to necessary conditions which a pair $(\mathcal M,M)$ satisfies if $(A,\mathfrak a)$ is an indecomposable (weak) coideal.

In Sections 4 and 5 the above mentioned general approach is applied to the problem of classification of
$\mathfrak G$-algebras and weak coideals of WHA's associated with a concrete class of fusion categories - Tambara-Yamagami categories $\mathcal T\mathcal Y(G,\chi,\tau)$ \cite{TY}. 

Recall that    simple objects of $\mathcal T\mathcal Y(G,\chi,\tau)$ are exactly the elements of a finite abelian group $G$ and one separate element $m$ satisfying the fusion rule $g\cdot h=gh,\ g\cdot m=m\cdot g=m,\ m^2=\underset{g\in G}\Sigma g,\ g^*=-g,\ m=m^*\ (g,h\in G)$. These categories are parameterized by non degenerate symmetric bicharacters $\chi:G\times G\to\mathbb C\backslash\{0\}$ and $\tau=\pm |G|^{-1/2}$. For any subset $K\subset G$, we shall denote $K^\perp:=\{g\in G|\chi(k,g)=1,\ \forall k\in K\}.$

The Hayashi's reconstruction theorem allows to construct a WHA  $\mathfrak G_{\mathcal T\mathcal Y}$ associated with $\mathcal T\mathcal Y(G,\chi,\tau)$ - see \cite{M}. We recall this construction in slightly different form in Subsection 4.1. Then, using the methods elaborated in \cite{MeMu}, we classify in Subsection 4.2 indecomposable module categories over representations of $\mathfrak G_{\mathcal T\mathcal Y}$, their autoequivalences and generators. Together with the above mentioned results this leads to the following classification theorem:
\begin{theorem} \label{TY} There are two types of isomorphism classes of indecomposable finite dimensional
$\mathfrak G_{\mathcal T\mathcal Y}$-$C^*$-algebras:

{\bf (i)} those parameterized by pairs $(K,\{m_\lambda\}^{orb})$, where $K<G$ and $\{m_\lambda\}^{orb}$ is the orbit of a nonzero collection $\{m_\lambda\in\mathbb Z_+|\lambda\in G/K\}$ under the action of the group of translations on $G/K$.

{\bf (ii)} those parameterized by pairs $(K, (\{m_\lambda\} ,\{m_\mu\})^{orb})$, where $K<G$ and $(\{m_\lambda\} ,\{m_\mu\})^{orb}$ is the orbit of a nonzero double collection $(\{m_\lambda\in\mathbb Z_+|\lambda\in G/K\},\{m_\mu\in\mathbb Z_+|\mu\in G/K^\perp\})$ under the action of:

a) the group of translations on $G/K\times G/K^\perp$ if $K\neq K^\perp$;

b) the semi-direct product $(G/K\times G/K)\underset{\sigma}\ltimes\mathbb Z_2$ generated by translations on $G/K\times G/K$ and the flip $\sigma:(\{m_\lambda\},\{m_\mu\})\leftrightarrow(\{m_\mu\},\{m_\lambda\})$ if $K= K^\perp$.
\end{theorem}

Finally, Section 5 is devoted to the classification of indecomposable (weak) coideals of $\mathfrak G_{\mathcal T\mathcal Y}$. Their classification is given by the following

\begin{theorem} \label{clas3} Isomorphism classes of indecomposable weak coideals of $\mathfrak G_{\mathcal T\mathcal Y}$ are parameterized by pairs $(K,(Z_0,Z_1)^{orb})$, where $K$ is a subgroup of $G$ and $(Z_0,Z_1)^{orb}$ is the orbit of a nonempty subset $(Z_0,Z_1)\subset G/K\times G/K^\perp$ such that either $|Z_0|\leq 1$ or $|Z_1|\leq 1$, under the action of:

a) the group of translations on $G/K\times G/K^\perp$ if $K\neq K^\perp$;

b) the semi-direct product $(G/K\times G/K)\underset{\sigma}\ltimes\mathbb Z_2$ generated by translations on $G/K\times G/K$ and the flip $\sigma:(Z_0,Z_1)\leftrightarrow (Z_1,Z_0)$ if $K=K^\perp$.

Given a subgroup $K<G$, the isomorphism classes containing coideals correspond exactly to the following orbits:

when $K\not= K^\perp$, to the four orbits $\{(\lambda,\emptyset) / \lambda \in G/K\}$,  $\{(\emptyset,\mu), / \mu \in G/K^\perp \}$, $\{(G/K,\mu) / \mu \in G/K^\perp \}$, $\{(\lambda,G/K^\perp) / \lambda \in G/K\}$,

when $K = K^\perp$, to the two orbits $\{(\lambda,\emptyset)\cup(\emptyset,\lambda),/ \lambda \in G/K \}$ and
$\{(G/K,\lambda) \cup (\lambda,G/K) / \lambda \in G/K\}$.

\end{theorem}

In fact, we give an explicit construction of representatives of all isomorphism classes of indecomposable finite dimensional
$\mathfrak G_{\mathcal T\mathcal Y}$-$C^*$-algebras and indecomposable (weak) coideals of $\mathfrak G_{\mathcal T\mathcal Y}$.

Our references are: to \cite{EGNO} for tensor categories, to \cite{Nes} for $C^*$-tensor categories and to \cite{NV} for weak Hopf algebras (finite quantum groupoids).
\end{section}
\begin{section}{Preliminaries}
\begin{subsection}{Weak Hopf $C^*$-algebras} A {\em weak Hopf  $C^*$-algebra} (WHA) $\mathfrak
G=(B,\Delta,S,\varepsilon)$ is a finite dimensional $C^*$-algebra $B$ with the
comultiplication $\Delta : B\to B\otimes B$, counit $\varepsilon : B\to\mathbb C$, and
antipode $S:B\to B$ such that $(B,\Delta, \varepsilon)$ is a coalgebra and the
following axioms hold for all $b,c,d\in B$ :
\begin{enumerate}
\item[(1)] $\Delta$ is a (not necessarily unital)  $*$-homomorphism :
$$
\Delta(bc) = \Delta(b)\Delta(c), \quad \Delta(b^*) =
\Delta(b)^*,
$$
\item[(2)] The unit and counit satisfy the identities (we use the Sweedler leg notation $\Delta(c)=c_1\otimes c_2,\
(\Delta\otimes id_B)\Delta(c)=c_1\otimes c_2\otimes c_3$ etc.):
\begin{eqnarray*}
\varepsilon(bc_1)\varepsilon(c_2d) &=&\varepsilon(bcd), \\
(\Delta(1)\otimes 1)(1\otimes \Delta(1)) &=& (\Delta\otimes id_B)\Delta(1),
\end{eqnarray*}
\item[(3)]
$S$ is an anti-algebra  and anti-coalgebra  map such that
\begin{eqnarray*}
m(id_B \otimes S)\Delta(b) &=& (\varepsilon\otimes id_B)(\Delta(1)(b\otimes 1)),\\
m(S\otimes id_B)\Delta(b) &=& (id_B \otimes \varepsilon)((1\otimes b)\Delta(1)),
\end{eqnarray*}
where $m$ denotes the multiplication.
\end{enumerate}
\medskip

The right hand sides of two last formulas are called {\em target}
and {\em source counital maps} $\varepsilon_t$ and $\varepsilon_s$,
respectively. Their images are unital $C^*$-subalgebras of $B$ called
{\em target} and {\em source counital subalgebras} $B_t$ and $B_s$,
respectively.

The dual vector space $\hat B$ has a natural structure of a weak Hopf
$C^*$-algebra $\hat{\mathfrak G}=(\hat B,\hat\Delta,\hat S,\hat\varepsilon)$
given by dualizing the structure operations of  $B$:
\begin{eqnarray*}
<\varphi\psi,\, b> &=& < \varphi\otimes\psi,\, \Delta(b)>, \\
<\hat\Delta(\varphi),\, b\otimes c> &=& <\varphi,\, bc>, \\
<\hat S(\varphi),\, b> &=& < \varphi,\, S(b)>, \\
< \phi^*,b> &=& \overline{ < \varphi,\, S(b)^*> },
\end{eqnarray*}
for all $b,c\in B$ and $\varphi,\psi\in \hat B$. The unit of $\hat B$ is
$\varepsilon$ and the counit is $1$.

The counital subalgebras commute elementwise, we have $S\circ\varepsilon_s =
\varepsilon_t\circ S$ and $S(B_t) =B_s$. We say that $B$ is {\em connected} if
$B_t \cap Z(B)= \mathbb{C}$ (where $Z(B)$ is the center of $B$), coconnected if
$B_t \cap B_s = \mathbb{C}$, and {\em biconnected} if both conditions are
satisfied.
\medskip

The antipode $S$ is unique, invertible, and satisfies $(S\circ *)^2 =id_B$. We
will only consider {\em regular} quantum groupoids, i.e., such that
$S^2|_{B_t}=id$. In this case, there exists a canonical positive element $H$ in
the center of $B_t$ such that $S^2$ is an inner automorphism implemented by $G=
HS(H)^{-1}$, i.e., $S^2(b) = GbG^{-1}$ for all $b\in B$. The element $G$ is
called the canonical group-like element of $B$, it satisfies the relation\
$\Delta(G) =(G\otimes G)\Delta(1)= \Delta(1)(G\otimes G)$.

There exists a unique positive functional $h$ on $B$,
called a {\em normalized Haar measure} such that
$$
(id_B\otimes h)\Delta = (\varepsilon_t\otimes h)\Delta,\quad
h\circ S =h,\quad h\circ \varepsilon_t = \varepsilon,\quad (id_B\otimes h)\Delta(1_B) = 1_B.
$$
We will dehote by $H_h$ the GNS Hilbert space generated by $B$ and $h$ and by $\Lambda_h:
B\to H_h$ the corresponding GNS map.
\end{subsection}
\begin{subsection}{Unitary representations and corepresentations of a weak Hopf $C^*$-algebra}

Let $\mathfrak G=(B,\Delta, S,\varepsilon)$ be a weak Hopf $C^*$-algebra. We denote by $\varepsilon_t, \varepsilon_s$
its target and source counital maps, by $B_t$ and $B_s$ its target and source subalgebras, respectively, and by $G$ its canonical group-like element. We also denote by $h$ the normalized Haar measure of $\mathfrak G$.

Any object of the category $URep(\mathfrak G)$ of unitary representations of $\mathfrak G$ is a left $B$-module of finite rank
such that the underlying vector space is a Hilbert space $H$ with a scalar product $<\cdot,\cdot>$ satisfying
$$
<b\cdot v,w>=<v,b^*\cdot w>,\quad\text{for all}\quad v,w\in H,\ b\in B.
$$
$URep(\mathfrak G)$ is a semisimple category whose morphisms are $B$-linear maps and simple objects are irreducible $B$-modules.
One defines the tensor product of two objects $H_1,H_2\in URep(\mathfrak G)$ as the Hilbert subspace $\Delta(1_B)\cdot(H_1\otimes H_2)$
of the usual tensor product together with the action of $B$ given by $\Delta$. Here we use the fact that $\Delta(1_B)$ is an orthogonal projection.

Tensor product of morphisms is the restriction of the usual tensor product of $B$-module morphisms. Let us note that any $H\in
URep(\mathfrak G)$ is automatically a $B_t$-bimodule via $z\cdot v\cdot t:=zS(t)\cdot v,\ \forall z,t\in B_t, v\in E$, and that
the above tensor product is in fact $\otimes_{B_t}$, moreover the $B_t$-bimodule structure for $H_1\otimes_{B_t} H_2$ is given by
$z\cdot \xi \cdot t =(z\otimes S(t))\cdot \xi,\ \forall z,t\in B_t, \xi \in H_1\otimes_{B_t} H_2$. The above tensor product is
associative, so the associativity isomorphisms are trivial. The unit object of $URep(\mathfrak G)$ is $B_t$ with the action of $B$
given by $b\cdot z:=\varepsilon _t(bz),\ \forall b\in B, z\in B_t$ and the scalar product $<z,t>=h(t^*z)$.

For any morphism $f:H_1\to H_2$, let $f^*:H_2\to H_1$ be the adjoint linear map: $<f(v),w>=<v,f^*(w)>,\ \forall v\in H_1, w\in H_2$.
Clearly, $f^*$ is $B$-linear, $f^{**}=f$, $(f\otimes_{B_t} g)^*=f^*\otimes_{B_t} g^*$, and $End(H)$ is a $C^*$-algebra, for any
object $H$. So $URep(\mathfrak G)$ is a finite $C^*$-multitensor category (${\bf 1}$ can be decomposable).

The conjugate object for any $H\in URep(\mathfrak G)$ is the dual vector space $\hat H$ naturally identified
($v\mapsto \overline v$) with the conjugate Hilbert space $\overline H$ with the action of $B$ defined by $b\cdot\overline v=\overline{G^{1/2}S(b)^*G^{-1/2}\cdot v}$, where $G$ is the canonical group-like element of $\mathfrak G$. Then the
rigidity morphisms defined by
\begin{equation} \label{rigid}
R_H(1_B)=\Sigma_i (G^{1/2}\cdot \overline e_i\otimes_{B_t}\cdot e_i),\ \overline R_H(1_B)=\Sigma_i (e_i\otimes_{B_t} G^{-1/2}\cdot\overline e_i),
\end{equation}
where $\{e_i\}_i$ is any orthogonal basis in $H$, satisfy all the needed properties - see \cite{BSz}, 3.6. Also,
it is known that the $B$-module $B_t$ is irreducible if and only if $B_s\cap Z(B)=\mathbb C 1_B$, i.e., if $\mathfrak G$ is
connected. So that, we have

\begin{proposition} \label{C*}
$URep(\mathfrak G)$ is a rigid finite $C^*$-multitensor category with trivial associativity
constraints. It is $C^*$-tensor if and only if $\mathfrak G$ is connected.
\end{proposition}

\begin{definition} \label{corepresentation}
A {\bf right unitary corepresentation} $U$ of $\mathfrak G$ on a Hilbert space $H_U$ is a partial isometry $U\in B(H_U)\otimes B$ such that:

(i) $(id\otimes\Delta)(U)=U_{12}U_{13}$.

(ii) $(id\otimes\varepsilon)(U)=id$.

If $U$ and $V$ are two right corepresentations on Hilbert spaces $H_U$ and $H_V$, respectively, a morphism between them is a
bounded linear map $T\in B(H_U,H_V)$ such that $(T \otimes 1_B)U= V(T\otimes 1_B)$. We denote by $UCorep(\mathfrak G)$ the
category whose objects are unitary corepresentations on {\bf finite dimensional} vector spaces with morphisms as above.
\end{definition}

If $\mathfrak G$ is coconnected (i.e., if $B_t\cap B_s=\mathbb C 1_B$), $UCorep(\mathfrak G)$ is a rigid $C^*$-tensor category with
trivial associativities isomorphic to $URep(\hat{\mathfrak G})$. Namely, any $H_U$ is a right $B$-comodule via $v\mapsto U(v\otimes 1_B)$,
therefore, automatically a $(B_s,B_s)$-bimodule. Then tensor product $U\otop V:=U_{13}V_{23}$ acts on $H_U\otimes_{B_s} H_V$, the unit
object $U_\varepsilon\in B(B_s)\otimes B$ is defined by $z\otimes b\mapsto \Delta(1_B)(1_B\otimes zb),
\ \forall z\in B_s, b\in B$, and the rigidity morphisms related to the conjugate $\overline U$ of an object $U$ which acts on the
conjugate Hilbert space $\overline H_U$ of $H_U$, are
\begin{equation} \label{rigid}
R_U(1_B)=\Sigma_i ({\hat G}^{1/2}\cdot \overline e_i\otimes_{B_s} e_i),\ \overline R_U(1_B)=\Sigma_i (e_i\otimes_{B_s}
{\hat G^{-1/2}}\cdot\overline e_i),
\end{equation}
where $\{e_i\}_i$ is any orthogonal basis in $H_U$. We denote by $\Omega$ an exhaustive set of representatives of the equivalence
classes of irreducibles in $UCorep(\mathfrak G)$.

Denote $H_{U^x}$ by $H^x$, then $U^x=\oplus_{i,j} m^x_{i,j}\otimes U^x_{i,j}$, where $m^x_{i,j}$ are the matrix units of $B(H^x)$
with respect to some orthogonal basis $\{e_i\}\in H^x$ and $U^x_{i,j}$ are the corresponding matrix coefficients of $U^x$.
Recall that $B=\oplus_{x\in\Omega} B_{U^x}$, where $B_{U^x}=Span(U^x_{i,j})$.
\end{subsection}
\begin{subsection}{ The Hayashi's fiber functor and reconstruction theorem.}

Let $\mathcal C$ be a rigid finite $C^*$-tensor category and $\Omega=Irr(\mathcal C)$ be an exhaustive set of representatives of equivalence
classes of its simple objects. Let $R$ be the $C^*$-algebra $R=\mathbb C^{\Omega}=\underset {x\in\Omega}\bigoplus\mathbb C p_x$, where $p_x=p_x^*$
are mutually orthogonal idempotents: $p_x p_y=\delta_{x,y}p_x$, for all $x,y\in \Omega$. Let us define a functor $\mathcal H$ from $\mathcal C$
to the category $Corr_f(R)$ of finite dimensional Hilbert $R$-bimodules (called $R$-correspondences) by:
$$
\mathcal H(x)=H^x=\underset {y,z\in\Omega}\bigoplus Hom(z,y\otimes x),\quad\text{for\ every}\ x\in \Omega,
$$
where $Hom(x,y)$ is the vector space of morphisms $x\to y$. The $R$-bimodule structure on $H^x$ is given by:
$$
p_y\cdot H^x\cdot p_z=Hom(z,y\otimes x),\quad\text{for\ all}\quad x,y,z\in\Omega.
$$
If $f\in Hom(x,y)$, then $\mathcal H(f):H^x\to H^y$ is defined by:
$$
\mathcal H(f)(g)=(id_z\otimes f)\cdot g,\quad\text{for\ any}\ z,t\in\Omega\ and\ g\in p_z\cdot H^x\cdot p_t.
$$
The tensor structure of $\mathcal H$ is a family of natural isomorphisms $\mathcal H_{x,y}: H^x\otimes H^y
\to H^x\underset {R}\otimes H^y$ defined by:
\begin{equation}\label{tensorhom}
\mathcal H_{x,y}(v\otimes w)=a_{z,x,y}\cdot(v\otimes id_y)\cdot w\in p_z\cdot H^{(x\otimes y)}\cdot p_s,
\end{equation}
for all $v\in p_z\cdot H^x\cdot p_t, w\in p_t\cdot H^y\cdot p_s, z,s,t\in\Omega$. Here $a_{z,x,y}$ are the
associativity isomorphisms of $\mathcal C$.

We define the scalar product on $H^x$ as follows. If $x,y,z \in \Omega$ and $f,g\in Hom(z,y\otimes x)$, then $g^*\in Hom(y\otimes x,z)$
and $g^*\cdot f\in End(z)=\mathbb C$, so one can put $<f,g>_x= g^*\cdot f $. The subspaces $Hom(z,y\otimes x)$ are declared to be orthogonal,
so $H^x\in Corr_f(R)$. Dually, $\overline H^x\in Corr_f(R)$ via $z_1\cdot{\overline v}\cdot z_2=\overline{z_2^*\cdot v\cdot z_1^*}$, for all $z_1,z_2\in R,v\in H^x$. Now one can check that $\mathcal H$ is a unitary tensor functor in the sense of \cite{Nes} 2.1.3.

\begin{theorem} \label{reconstruction} (a $C^*$-version of the Hayashi's theorem -see \cite{Ha}, \cite{Pf})

Let $\mathcal C$ be a rigid finite $C^*$-tensor category, $\Omega=Irr(\mathcal C)$ and $\mathcal H : \mathcal C \to Corr_f(R)$ be the Hayashi's functor, where $R=\mathbb C^{|\Omega}|$. Then the vector space
\begin{equation} \label{algebra}
B =  \underset {x \in \Omega} \bigoplus H_x\otimes\overline  H_x,
\end{equation}
has a regular biconnected weak Hopf $C^*$-algebra structure $\mathfrak G$ such that $\mathcal C\cong UCorep(\mathfrak G)$ as rigid $C^*$-tensor categories.
\end{theorem}
Explicitly, if $v,w\in H^x, g,h\in H^y$ and $\{e^{x}_j\}$ is an orthonormal basis in $H^x$, for all $x,y\in\Omega$, then:
\begin{equation} \label{product}
(w\otimes\overline v)_x\cdot (g\otimes\overline h)_y =  ({\mathcal H_{x,y}(w\otimes g)}\otimes \overline{\mathcal H_{x,y}
(v\otimes h)})_{x\otimes y}\in H^{(x\otimes y)}\otimes\overline{H^{(x\otimes y)}}
\end{equation}
\begin{equation} \label{coproduct}
\Delta(w\otimes\overline v) =  \underset {j}\bigoplus (w\otimes\overline {e^{x}_j})_x\otimes(e^{x}_j\otimes\overline  v)_x,
\end{equation}
\begin{equation} \label{counit}
\varepsilon (w\otimes\overline v) =<w,v>_x.
\end{equation}
Now define an antipode and an involution. Consider the natural isomorphisms $\Phi_x:H^x\to\overline H^{x^*}$
and $\Psi_x:\overline H^x\to H^{x^*}$, where $x^*$ is the dual of $x\in\Omega$:
\begin{equation} \label{starinv}
\Phi_x(v)=(id_y\otimes \overline{R_x}^*)\cdot a_{y,x,x^*}\cdot(v\otimes id_{x^*}), \Psi_x(\overline v)=
(\overline v\otimes id_{x^*})\cdot a^{-1}_{y,x,x^*}\cdot(id_y\otimes\overline{R_x}),
\end{equation}
where $x,y,z\in\Omega$, we identify $y$ with $y\otimes 1$, $v\in {p_y}\cdot H^x\cdot {p_z}$,
$\overline{R_x}$ and $a_{y,x,x^*}$ are, respectively, the rigidity morphisms and the associativities in $\mathcal C$. Then:
\begin{equation} \label{antipode}
S(w\otimes \overline v)=\Psi_x(\overline v)\otimes\Phi_x(w),
\end{equation}
\begin{equation} \label{involution}
(w\otimes \overline v)^*=w^\natural\otimes{\overline v}^\flat,\ \text{where}\ w^\natural=\Psi_x(\overline w),\ {\overline v}^\flat=\Phi_x(v).
\end{equation}
Any $H^x$ is a right $B$-comodule via
$$
\mathfrak a_x(v)= \underset {j}\Sigma e^x_j \otimes \overline{e^x_j}\otimes v ,\quad \text{where}\quad v\in H_x,
$$
one checks that it is unitary which gives the equivalence $\mathcal C\cong UCorep(\mathfrak G)$.

The algebra of the dual quantum groupoid $\hat{\mathfrak G}$ is
\begin{equation} \label{dualalg}
\hat B =  \underset {x\in\Omega}\bigoplus B(H_x),
\end{equation}
the duality is given, for all $x\in\Omega, A\in B(H_x),v,w\in H_x$ by:
$$
<A,w\otimes\overline v>=<Aw,v>_x.
$$
$\hat B$ is clearly a $C^*$-algebra with the obvious matrix product and involution,

\begin{notations}\label{circ}
 For all $x,y \in \Omega$ and all $v \in H^x$, $w \in H^y$, we denote:
 $$ v\circ w = \mathcal H_{x,y}(v \otimes_R w)$$
\end{notations}

\begin{remark}
\label{hzero}
Let $0$ be the unit element of $\mathcal C$,  and  $H^0 := \underset{x \in \Omega} \oplus Hom(x,x)$, then using (\ref{tensorhom}) and (\ref{starinv}) it is easy to check that $(H^0,\circ,\sharp)$ is a commutative $C^*$-algebra
and if, for all $x \in \Omega$, $v^0_x$ is a normalized vector in $ Hom(x,x)$, then $(v^0_x)_{x \in \Omega}$ is a basis
of mutually orthogonal projections in $H^0$.
\end{remark}

\begin{remark} \label{rec}
Let $\mathcal C$ be a rigid finite $C^*$-tensor category and $\mathcal F : \mathcal C \to Corr_f(R)$ be a unitary tensor functor,
where $R$ is a finite dimensional unital $C^*$-algebra. Then there exists \cite{Sz} a regular biconnected finite quantum groupoid $\mathfrak G$ with $B_t\cong B_s\cong R$ such that $\mathcal C\cong UCorep(\mathfrak G)$ as $C^*$-tensor categories. For any fixed $\mathcal C$, the set of $C^*$-algebras $R$ for which the above mentioned functor $\mathcal F$ exists, contains at least $R=\mathbb C^{|\Omega|}$ (where $\Omega=Irr(\mathcal C)$), then $\mathcal F=\mathcal H$. In general, this set contains several elements, and the corresponding WHAs are called Morita equivalent.

In particular, if the above set of functors contains a {\bf fiber} functor $\mathcal F : \mathcal C \to Hilb_f$, i.e., $R=\mathbb C$, the corresponding quantum groupoids are Morita equivalent to a usual $C^*$-Hopf algebra.
\end{remark}
\end{subsection}


\begin{subsection}{Coactions.}
\begin{definition}\label{sieste}
A right coaction of a WHA $\mathfrak G$ on a unital $*$-algebra $A$, is a $*$-homomorphism
$\mathfrak a: A \to A \otimes B$  such that:


1) $(\mathfrak a \otimes i)\mathfrak a= (id_A\otimes \Delta)\mathfrak a$.

2) $(id_A\otimes\varepsilon )\mathfrak a=id_A.$

3) $\mathfrak a(1_A) \in A \otimes B_t$.

One also says that $(A,\mathfrak a)$ is a $\mathfrak G$-$*$-algebra.
\end{definition}
If $A$ is a $C^*$-algebra, then $\mathfrak a$ is automatically continuous, even an isometry.

There are $*$-homomorphism $\alpha:B_s\to A$ and $*$-antihomomorphism $\beta:B_s\to A$ with commuting images defined by
$\alpha(x)\beta(y):=(id_A\otimes\varepsilon)[(1_A\otimes x)\mathfrak a(1_A)(1_A\otimes y)]$, for all $x,y\in B_s$.
We also have $\mathfrak a(1_A)=(\alpha\otimes id_B)\Delta(1_B)$,
\begin{equation} \label{compatible}
\mathfrak a(\alpha(x)a\beta(y))=(1_A\otimes x)\mathfrak a(a)(1_A\otimes y),
\end{equation}
and
\begin{equation} \label{compatible'}
(\alpha(x)\otimes 1_B)\mathfrak a(a)(\beta(y)\otimes 1_B)=(1_A\otimes S(x))\mathfrak a(a)(1_A\otimes S(y)).
\end{equation}
The set $A^{\mathfrak a} = \{a \in A | \mathfrak a(a) = \mathfrak a(1_A) (a \otimes 1_B) \}$ is a unital $*$-subalgebra
of $A$ (it is a unital $C^*$-subalgebra of $A$ when $A$ is a $C^*$-algebra) commuting pointwise with $\alpha(B_s)$.
A coaction $\mathfrak a$ is called {\it ergodic} if $A^{\mathfrak a} = \mathbb C 1_A$.

\begin{definition}\label{indecomp}
A $\mathfrak G-C^*$-algebra $(A,\mathfrak a)$ is said to be indecomposable if it cannot be presented as
a direct sum of two $\mathfrak G-C^*$-algebras.
\end{definition}
It is easy to see that $(A,\mathfrak a)$ is indecomposable if and only if $Z(A)\cap A^{\mathfrak a}=\mathbb C 1_A$.
Clearly, any ergodic $\mathfrak G-C^*$-algebra is indecomposable.

For any $(U,H_U)\in UCorep(\mathfrak G)$, we define the {\it spectral subspace} of $A$ corresponding to $(U,H_U)$ by
$$
A_U:=\{a\in A|\mathfrak a(a)\in \mathfrak a(1_A)(A\otimes B_U)\}.
$$
Let us recall the properties of the spectral subspaces:

(i) All $A_U$ are closed.

(ii) $A=\oplus_{x\in \Omega} A_{U^x}$.

(iii) $A_{U^x}A_{U^y}\subset \oplus_z A_{U^z},$ where $z$ runs over the set of all irreducible direct
summands of $U^x\otop U^y$.

(iv) $\mathfrak a(A_U)\subset \mathfrak a(1_A)(A_U\otimes B_U)$ and $A_{\overline U}=(A_U)^*$.

(v) $A_\varepsilon $ is a unital $C^*$-algebra.
\end{subsection}
\begin{subsection}{Categorical duality.}
Let us recall the main result of \cite{VV1}:
\begin{theorem} \label{main}
Given a regular coconnected WHA $\mathfrak G$, the following two categories are equivalent:

(i) The category of unital $\mathfrak G$-$C^*$-algebras with unital $\mathfrak G$-equivariant $*$-homomorphisms as morphisms.

(ii) The category of pairs $(\mathcal M, M)$, where $\mathcal M$ is a left module $C^*$-category with trivial module associativities over the $C^*$-tensor category $UCorep(\mathfrak G)$ and $M$ is a generator in $\mathcal M$, with equivalence classes of unitary module functors respecting the prescribed generators as morphisms.
\end{theorem}
In particular, given a unital $\mathfrak G$-$C^*$-algebra $A$, one constructs the $C^*$-category $\mathcal M=\mathcal D_A$ of
finitely generated right Hilbert $A$-modules which are equivariant, that is, equipped with a compatible right coaction \cite{BS1}.
Any its object is automatically a $(B_s,A)$-bimodule, and the bifunctor $U\boxtimes X:= H_U\otimes_{B_s} X\in \mathcal D_A$,
for all $U\in UCorep(\mathfrak G)$ and $X\in\mathcal D_A$, turns $\mathcal D_A$ into a left module $C^*$-category over
$UCorep(\mathfrak G)$ with generator $A$ and trivial associativities.

Vice versa, if a pair $(\mathcal M, M)$ is given, the construction of a $\mathfrak G$-$C^*$-algebra $(A,\mathfrak a)$
contains the following steps. First, denote by $R$ the unital $C^*$-algebra $End(M)$ and consider the functor $F:\mathcal C\to Corr(R)$
defined on the objects by $F(U)=Hom_{\mathcal M} (M, U\boxtimes M)\ \forall U\in\mathcal C$. Here $X=F(U)$ is a right $R$-module via the
composition of morphisms, a left $R$-module via $rX=(id\otimes r)X$, the $R$-valued inner product is given by $<X,Y>=X^*Y$, the action of
$F$ on morphisms is defined by $F(T)X=(T\otimes id)X$. The weak tensor structure of $F$ (in the sense of \cite{Nes1}) is given by $J_{X,Y}
(X\otimes Y)=(id\otimes Y)X$, for all $X\in F(U), Y\in F(V), U,V\in UCorep(\mathfrak G)$.

Then consider two vector spaces:
\begin{equation} \label{bspace}
A=\underset{x\in\Omega}\bigoplus A_{U^x}:=\underset{x \in\Omega}\bigoplus( F(U^x)\otimes\overline{H^x})
\end{equation}
and
\begin{equation} \label{recalgebra2}
\tilde{A}=\underset{U\in \|UCorep(\mathfrak G)\|}\bigoplus A_{U}:=\underset{U\in \|UCorep(\mathfrak G)\|}
\bigoplus(\overline{F(U)}\otimes H_U),
\end{equation}
where $F(U)=\underset{i}\bigoplus F(U_i)$ corresponds to the decomposition $U=\bigoplus U_i$ into irreducibles, and
$\|UCorep(\mathfrak G)\|$ is an exhaustive set of representatives of the equivalence classes of objects in
$UCorep(\mathfrak G)$ (these classes constitute a countable set). $\tilde A$ is a unital associative algebra with
the product
$$
(X\otimes\overline\xi)(Y\otimes\overline\eta)=(id\otimes Y)X\otimes(\overline\xi\otimes_{B_s}\overline\eta),
\ \forall (X\otimes\overline\xi)\in A_U, (Y\otimes \overline\eta)\in A_V,
$$
and the unit
$$
1_{\tilde A}=id_M\otimes\overline{1_B}.
$$
Note that $(id\otimes Y)X=J_{X,Y}(X\otimes Y)\in F(U\otop V)$. Then, for any $U\in UCorep(G)$, choose isometries
$w_i:H_i\to H_U$ defining the decomposition of $U$ into irreducibles, and construct the projection $p:\tilde A\to A$ by
\begin{equation} \label{p}
p(X\otimes\overline\xi)=\Sigma_i(F(w_i^*)X\otimes\overline{w_i^*\xi}),\ \forall (X\otimes\overline\xi)\in A_U,
\end{equation}
which does not depend on the choice of $w_i$. Then $A$ is a unital $*$-algebra with the product $x\cdot y:=
p(xy)$, for all $x,y\in A$ and the involution  $x^*:=p(x^\bullet)$, where $(X\otimes\overline\xi)^\bullet:=
(id\otimes X^*) F(\overline R_U)\otimes\overline{\hat G^{1/2}\xi}$, for all $\xi\in H_U, X\in F(U), U\in UCorep(\mathfrak G)$.
Here $\overline R_U$ is the rigidity morphism from (\ref{rigid}). Finally, the map
\begin{equation} \label{coact}
\mathfrak a(X\otimes\overline\xi_i)=X\otimes\Sigma_{j}(\overline\xi_j\otimes U^x_{j,i}),
\end{equation}
where $\{\xi_i\}$ is an orthogonal basis in $H^x$ and $(U^x_{i,j}$ are the matrix elements of $U^x$ in this basis, is
a right coaction of $\mathfrak G$ on $A$. Moreover, $A$ admits a unique $C^*$-completion $\overline A$ such that $\mathfrak a$
extends to a continuous coaction of $\mathfrak G$ on it.
\begin{remark} \label{decomp}
1) We say that a $UCorep(\mathfrak G)$-module category is indecomposable if it is not equivalent to a direct sum of two nontrivial
$UCorep(\mathfrak G)$-module subcategories. Theorem \ref{main} implies that a $\mathfrak G-C^*$-algebra $(A,\mathfrak a)$ is
indecomposable if and only if the $UCorep(\mathfrak G)$-module category $\mathcal M$ is indecomposable.

2) Let $I$ be a unital right coideal $*$-subalgebra of $B$. Then $I_{U^x}=I\cap B_{U^x}$ and $F(U^x)$ can be
identified with a Hilbert subspace of $H^x\ (\forall x\in\Omega)$ and the coaction is the restriction of $\Delta$.
\end{remark}
\begin{example} \label{regular}
The  $C^*$-algebra $B$ with coproduct $\Delta$ viewed as $\mathfrak G$-$C^*$-algebra, corresponds to the $UCorep(\mathfrak G)$-module $C^*$-category
$\mathcal Corr_f(B_s)$ with generator $M=B_s$: for any element $U\in UCorep(\mathfrak G)$ and $N\in Corr_f(B_s)$, one defines $U\boxtimes N:=
F(U)\otimes_{B_s} N$, where the functor $F: UCorep(\mathfrak G)\to Corr_f(B_s)\ (F(U)=H_U)$ is the forgetful functor. Indeed, identifying
$\mathcal M(B_s,H_U)$ with $H_U$, we get an isomorphism of the algebra $\tilde A$ constructed from the pair $(\mathcal M,M)$ onto
$\tilde B=\underset{U}\bigoplus(H_U\otimes\overline H_U)$ and then an isomorphism $A\cong B=\underset{x\in\hat G}\bigoplus (H_x\otimes\overline H_x)$
such that $p:\tilde A\to A$ turns into the map $\tilde B\to B$ sending $\xi\otimes\overline\eta\in H_U\otimes\overline H_U$ into the matrix coefficient $U_{\xi,\eta}$.
\end{example}
\end{subsection}
\end{section}
\begin{section} {Classifying Indecomposable Weak Coideals}

If $dim(A)<\infty$, we have the following remarks.

\begin{remark} \label{semisimple} If $(A,\mathfrak a)$ is a finite dimensional $\mathfrak G-C^*$-algebra, then
$\mathcal M=\mathcal D_A$ is a semisimple $C^*$-category. Indeed, $dim(Hom_{\mathcal M}(\mathcal E,\mathcal E))<\infty$,
for any $\mathcal E\in D_A$ which is finitely generated. Then the proof of \cite{CY}, Proposition 3.9 applies. As $A$
is a generator of $\mathcal M$, the set $\{M_\lambda|\lambda\in\Lambda\}$ of its (classes of) simple objects is
finite and we have the corresponding fusion rule
\begin{equation} \label{fusion}
U_x\boxtimes M_\lambda=\underset{\mu}\Sigma n^\mu_{x,\lambda} M_\mu,\ \text{where}\ x\in\Omega, n^\mu_{x,\lambda}=dim
(Hom_{\mathcal M}(U_x\boxtimes M_\lambda,M_\mu))\in \mathbb Z_+.
\end{equation}
The associativity and the unit object conditions mean, respectively, that
\begin{equation} \label{assocfusion}
\underset{z\in\Omega}\Sigma c^z_{x,y}n^\rho_{z,\lambda}= \underset{\mu\in\Lambda}\Sigma n^\rho_{x,\mu}n^\mu_{y,\lambda},
\ \text{and}\ n^\rho_{1,\lambda}=\delta_{\rho,\lambda},\quad\forall x,y\in\Omega, \rho,\lambda\in\Lambda,
\end{equation}
where $c^z_{x,y}$ are the fusion coefficients of $\mathcal C=UCorep(\mathfrak G)$. Proposition 7.1.6 of \cite{EGNO}
gives $n^\mu_{x,\lambda}=n^\lambda_{x^*,\mu}$, for all $\lambda,\mu\in\Lambda, x\in\Omega$.
\end{remark}

\begin{remark} \label{z_+-modules}
If $A$ is a coideal of $B$, then, due to \cite{VV}, Theorem 1.1, there is an inclusion $j:\mathcal M\mapsto \mathcal C$
such that
\begin{equation} \label{ucoid}
j(M)=\underset{x\in\Omega}\oplus U^x,
\end{equation}
where $\mathcal M$ is the left $\mathcal C$-module category with generator $M$ coming from $(A,\Delta|_A)$ and
$\mathcal C=UCorep(\mathfrak G)$ is viewed as a $\mathcal C$-module category with generator the
$\underset{x\in\Omega}\oplus U^x$.

If $\Lambda$ is the set of irreducibles of $\mathcal M$ (we denote them by $M_\lambda$), we can write $j(M_\lambda)=\underset{x\in\Omega}\Sigma a_{\lambda,x} U^x$, for all $\lambda\in\Lambda$,
where $a_{\lambda,x}\in\mathbb Z_+$.
Writing $M=\underset{\lambda\in\Lambda}\Sigma m_\lambda M_\lambda$, we must have:
\begin{equation} \label{mlambda}
\underset{\lambda\in\Lambda}\Sigma m_\lambda a_{\lambda,x}=1,\quad\text{for\ all}\quad x\in\Omega.
\end{equation}
\end{remark}

Recall that due to the reconstruction theorem for $\mathfrak G$, any $H^x (x\in\Omega)$ is the direct sum of 1-dimensional subspaces $Hom(z,y\otimes x)$, where $y,z\in\Omega$ are such that $z\subset (y\otimes x)$. In particular, $H^0=\underset{z\in\Omega}\oplus
Hom(z,z)$ (where $0$ denotes the trivial corepresentation of $\mathfrak G$); we will denote by $v^0_z$ a norm one vector generating
$Hom(z,z)$ viewed as a subspace of $H^0$.

The following lemma allows to select weak coideals of $B$ from   all $\mathfrak G-C^*$-algebras.

\begin{lemma} \label{wcoid} Let us fix a $UCorep(\mathfrak G)$-module category $\mathcal M$ and a generator $M$ in it, and
let $(A,\mathfrak a)$ be a $\mathfrak G$-algebra constructed from this data using the weak tensor functor $(F,J_{U,V})$. Then:

a) $(A,\mathfrak a)$ is a weak coideal of $B$ if and only if each $F(U^x)$ can be identified with a subspace $X^x\subset H^x$
such that the map $\zeta\mapsto\zeta^\natural=\Psi_x(\overline{\zeta})$ sends $X^x$ onto $X^{\overline x}\cong F(\overline{U^x})$
and $J_{U^x,U^y}=\mathcal H_{x,y}$, for all $x,y\in\Omega$.

b) $X^0$ is a $C^*$-subalgebra of $H^0$. The unit of $X^0$ is $v^0_\Gamma:=\underset{x\in\Gamma}\oplus v^0_x$, where $\Gamma\subset\Omega$
is some nonempty subset. $A=\underset{x\in\Omega}\oplus(X^x\otimes\overline H^x)$ is a coideal if and only if $\Gamma=\Omega$.

c) A weak coideal $A=\underset{x\in\Omega}\oplus(X^x\otimes\overline H^x)$ is decomposable if and only if
$Z(A)$ contains an element of the form $p=v^0_{\Gamma_0}\otimes v^0_\Omega$, where $\Gamma_0$ is a proper nonempty subset
of $\Gamma$.

d) For any two identifications, $F(U^x)\cong X^x$ and $F(U^x)\cong \tilde X^x, \forall x\in\Omega$, satisfying the above
mentioned conditions, the corresponding weak coideals $\underset{x\in\Omega}\oplus(X^x\otimes\overline H^x)$ and
$\underset{x\in\Omega}\oplus(\tilde X^x\otimes\overline H^x)$ are isomorphic as $\mathfrak G$-$C^*$-algebras.

\end{lemma}

{\bf Proof.} a) If $(A,\mathfrak a)$ is a weak coideal of $B$, then $A_U\subset B_U$, for any $U\in UCorep(\mathfrak G)$.
Indeed, by \cite{VV1}, Proposition 3.17 $A_U=\{a\in A|\Delta(a)\in \Delta(1_A)(A\otimes B_U)\}$, but
$\Delta(1_A)=\Delta(1_B)(1_A\otimes 1_B)$, hence $\Delta(a)\in \Delta(1_B)(A\otimes B_U)\subset \Delta(1_B)(B\otimes B_U)$, so that $A_U\subset B_U$. It follows from \cite{VV1}, Theorem 4.12 and Theorem 2.21, respectively, that $A_{U^x}\cong F(U^x)\otimes
\overline{H^x}$ and $B_{U^x}\cong H^x\otimes\overline{H^x}$, so the above inclusions mean that $F(U^x)\subset H^x$, for all $x\in\Omega$.

The multiplication in $A$ is the restriction of that in $B$, therefore, comparing the formulas (15) and \cite{VV1}, (16) and using
the relation $J_{U^x,U^y}(X\otimes_{B_s} Y)=(id\otimes Y)X (\forall X\in F(U^x), Y\in F(U^y))$, we have $J_{U^x,U^y}=\mathcal H_{x,y}$.

The involution in $A$ sends $X\otimes\overline\eta$ onto $X^*\otimes(\overline\eta)^\flat$ (see Subsection 2.4) and is the restriction
of that in $B$, the last one is defined by $\zeta\otimes\overline\eta\mapsto\zeta^\natural\otimes(\overline\eta)^\flat\ \forall \zeta,\eta
\in H^x, x\in\Omega$. Then for $X\in F(U^x)\subset H^x$ we have $X^*=X^\natural$.

Conversely, suppose that $F(U^x)\subset H^x$ and $J_{U^x,U^y}=\mathcal H_{x,y}$, for all $x,y\in\Omega$. It follows from the argument
above that the multiplication in $A$ is the restriction of that in $B$. Next, compare the formulas \cite{VV1}, (29) for $\mathfrak a$
and \cite{VV1}, (14) for $\Delta$. Since $B_{U^x}=H^x\otimes\overline H^x$, for any $U^x$ - see \cite{VV1}, (12), the matrix coefficient $U^x_{\zeta,\eta}$ with respect to a basis $\{\zeta^x\}$ of $H^x$ can be identified with $\zeta^x\otimes \eta^x$, for all $x\in\Omega$.
Now it is clear that $\mathfrak a$ is the restriction of $\Delta$. Finally, putting $X^*=X^\natural$ for any $X\in F(U^x)$ and using the
fact that $F(U^x)^\natural=F(\overline{U^x})$, one checks that $(A,\mathfrak a)$ is a coideal of $B$.

b) By Remark \ref{hzero}, $H^0=\underset{x\in\Omega}\oplus\mathbb C v^0_x$ is a commutative unital $C^*$-algebra, $v^0_x\ (x\in\Omega)$ are mutually orthogonal projections, and if $A=\underset{x\in\Omega}\oplus(X^x\otimes\overline H^x)$ is a weak coideal of $B$, then $X^0$ is a $C^*$-subalgebra of $H^0$. Its spectral mutually orthogonal projectors are $v^0_{\Gamma_i}$, where $\Gamma_i\subset\Omega\ (i=1,...,k_0=dim(X^0))$ are disjoint subsets of $\Omega$, the unit of $X^0$, i.e., the image
of $F(id_M)$, is $v^0_\Gamma$, where $\Gamma=\sqcup^{k_0}_{i=1}\Gamma_i$.
As $1_A=v^0_\Gamma\otimes\overline v^0_\Omega$ and $1_B=v^0_\Omega\otimes\overline v^0_\Omega$, $A$ is a coideal if and
only if $\Gamma=\Omega$.

c) One checks that $B_t=H^0\otimes \overline v^0_\Omega$ and that any nontrivial orthoprojector $p\in [Z(A)\cap B_t]$ gives a
decomposition $A=pA\oplus(1-p)A$ into the direct sum of two weak coideals of $B$. As $1_A=v^0_\Gamma\otimes\overline v^0_\Omega$,
$p$ must be of the form $v^0_{\Gamma_0}\otimes\overline v^0_\Omega$, where $\Gamma_0$ is a proper nonempty subset of $\Gamma$.

d) The two $\mathfrak G$-$C^*$-algebras are isomorphic because they correspond to the same couple $(\mathcal M,M)$.
\hfill $\square$
\begin{corollary} \label{multiplicity}
It follows from the definition of the functor $F$ that $X^0=F(U^0)=End_{\mathcal M}(M)$. This finite dimensional
$C^*$-algebra is commutative due to the statement b) which is only possible if $m_\lambda\in\{0,1\}$ for all $\lambda\in\Lambda$.
\end{corollary}
\end{section}


\begin{section}{Weak Hopf $C^*$-Algebras related to Tambara-Yamagami categories}
\label{sTY}

\subsection{ Tambara-Yamagami categories}

These categories denoted by $\mathcal T\mathcal Y(G,\chi,\tau)$ ($G$ is a finite group; we consider them only over $\mathbb C$) are $\mathbb Z_2$-graded fusion categories whose $0$-component is $Vec_G$ - the category of finite dimensional $G$-graded vector spaces with trivial associativities (its simple objects are $g\in G$) and $1$-component is generated by single simple object $m$. The Grothendieck ring of $\mathcal T\mathcal Y(G,\chi,\tau)$ is isomorphic to the $\mathbb Z_2$-graded fusion ring $\mathcal T\mathcal Y_G=\mathbb Z G\oplus\mathbb Z\{m\}$ such that $g\cdot m=m\cdot g=m,\ m^2=\underset{g\in G}
\Sigma g,\ m=m^*$. These categories exist if and only if $G$ is abelian, they are parameterized by non degenerate symmetric bicharacters
$\chi:G\times G\to\mathbb C\backslash\{0\}$ and $\tau=\pm |G|^{-1/2}$ - see \cite{TY}, \cite{EGNO}, Example 4.10.5. The associativities
$\phi(U,V,W):(U\otimes V)\otimes W\to U\otimes(V\otimes W)$ are
$$
\phi(g,h,k)=id_{g+h+k},\quad \phi(g,h,m)=id_{m},\quad \phi(m,g,h)=id_{m},
$$
$$
\phi(g,m,h)=\chi(g,h)id_{m},\quad \phi(g,m,m)=\underset{h\in G}\oplus id_{h},\quad \phi(m,m,g)=\underset{h\in G}\oplus id_{h},
$$
$$
\phi(m,g,m)=\underset{h\in G}\oplus \chi(g,h)id_{h},\quad \phi(m,m,m)=(\tau\chi(g,h)^{-1} id_{m})_{g,h},
$$
where $g,h,k\in G$. The unit isomorphisms are trivial. $\mathcal T\mathcal Y(G,\chi,\tau)$ becomes a $C^*$-tensor category
when $\chi:G\times G\to T=\{z\in\mathbb C||z|=1\}$, from now on we assume that this is the case. The dual objects are: $g^*=-g$,
for all $g\in G$, and $m^*=m$. The rigidity morphisms are defined by $R_g:0\overset{id_0}\to g^*\otimes g$, $\overline{R_g}:
0\overset{id_0}\to g\otimes g^*$, $R_m=\tau |G|^{1/2}\iota$, and $\overline{R_m}=|G|^{1/2}\iota$, where $\iota:0\to
m\otimes m$ is the inclusion. Then $dim_q(g)=1$, for all $g\in G$, and $dim_q(m)=\sqrt{|G|}$.

Let us apply Theorem \ref{reconstruction} to the category $\mathcal T\mathcal Y(G,\chi,\tau)$ in order to construct
a biconnected regular WHA $\mathfrak G_{\mathcal T\mathcal Y}=(B,\Delta,S,\varepsilon)$ with $UCorep(\mathfrak G_{\mathcal T
\mathcal Y})\cong\mathcal T\mathcal Y(G,\chi,\tau)$. The Hayashi's functor $\mathcal H:\mathcal T\mathcal Y(G,\chi,\tau)
\to Corr_f(R)$, where $C^*$-algebra $R:=End(\underset{x\in\Omega}\oplus x)\cong \mathbb C^{|G|+1}$, was constructed in
\cite{M}. Denoting $\Omega_g=\Omega:=G\sqcup\{m\}$ and $\Omega_m:=G\sqcup\overline G$, where $g\in G$ and
$\overline G$ is the second copy of $G$, one easily computes that $H^g\cong\mathbb C^{|G|+1}$, for all $g\in G$ and
$H^m:\cong\mathbb C^{2|G|}$.

Let us fix a basis $\{v^x_y\}(y\in\Omega_x)$ in each $H^x\ (x\in\Omega)$ choosing a norm one vector in every 1-dimensional vector subspace: $v^g_h\in Hom(h,(h-g)\otimes g)$, $v^g_m\in Hom(m,m\otimes g)$, $v^m_g\in Hom(m,g\otimes m)$, and $v^m_{\overline g}\in Hom(g,m\otimes m)$, where $g\in G$.

\begin{lemma}
\label{circsharp}
Using notations (\ref{involution}) and \ref{circ}, for all $g,h,k\in G$, $x \in \Omega$, one has:
$$
v^g_k\circ v^h_{x}= \delta_{x,h+k}v^{g+h}_{h+k},\ v^g_m\circ v^h_x= \delta_{x,m}v^{g+h}_m,
$$
$$
v^m_k\circ v^g_x=\delta_{x,m}\chi(g,k)v^m_k,\ v^m_{\overline k}\circ v^g_{x}=
\delta_{x,g+k}v^m_{\overline{g+k}},
$$
$$
v^g_x\circ v^m_{\overline k}=\delta_{x,m}\chi(g,k) v^m_{\overline k},\ v^g_x\circ v^m_k
= \delta_{x,k}v^m_{k-g},
$$
$$
v^m_h\circ v^m_{\overline k}= v^{k-h}_k,
\ v^m_{\overline h}\circ v^m_k=\delta_{h,k}\tau \underset{g\in G}\Sigma \chi(g,h)^{-1}v^g_m.
$$
$$
(v^g_k)^\sharp = v^{-g}_{k-g}, (v^g_m)^\sharp = v^{-g}_m, (v^m_g)^\sharp = |G|^{1/2}v^m_{\overline g}, (v^m_{\overline g})^\sharp = \tau^{-1}|G|^{1/2}v^m_g
$$
\end{lemma}
\begin{dm} For equations related to  product $\circ$,  these are computations made in \cite{M} 2.1.5, where $\mathcal H_{x,y}$ must be replaced by $\mathcal F_{x,y}^{-1}$. Moreover, in the case of $\mathcal T \mathcal Y(G,\chi,\tau)$, the formulas of (\ref{starinv}) imply that the isomorphisms $\Phi_x:H^x\to\overline{H^{x^*}}$ and $\Psi_x:\overline{H^x}\to H^{x^*}\ (x\in\Omega)$ of Theorem \ref{reconstruction}
are  given, for all $g,h \in G$,  by:
$$
\Phi_g(v^g_h)=\overline{v^{-g}_{h-g}},\ \Phi_g(v^g_m)=\overline{v^{-g}_m},\ \Phi_m(v^m_g)=|G|^{-1/2}\overline{v^m_{\overline g}},\
\Phi_m(v^m_{\overline g})=\tau|G|^{-1/2}\overline{v^m_g},
$$
$$
\Psi_g(\overline{v^g_h})=v^{-g}_{h-g},\ \Psi_g(\overline{v^g_m})=v^{-g}_m,\ \Psi_m(\overline{v^m_g})=|G|^{1/2}v^m_{\overline g},\
\Psi_m(\overline{v^m_{\overline g}})=\tau^{-1}|G|^{1/2}v^m_g,
$$
which implies, by (\ref{involution}) the formulas for involution $\sharp$.
\hfill\end{dm}


Now the whole structure of a WHA $\mathfrak G_{\mathcal T\mathcal Y}$ is given by formulas (\ref{algebra}), (\ref{product}),
(\ref{coproduct}), and (\ref{counit}). It was shown in \cite{M} that this WHA is isomorphic to its dual whose $C^*$-algebra
$\hat B\cong\underset{x\in\Omega}\oplus B(H^x)$. This implies that $URep(\mathfrak G_{\mathcal T\mathcal Y})\cong
UCorep(\mathfrak G_{\mathcal T\mathcal Y})$.

The isomorphisms $\Phi_x:H^x\to\overline{H^{x^*}}$ and $\Psi_x:\overline{H^x}\to H^{x^*}\ (x\in\Omega)$ of Theorem \ref{reconstruction}
are now given by:
$$
\Phi_g(v^g_h)=\overline{v^{-g}_{h-g}},\ \Phi_g(v^g_m)=\overline{v^{-g}_m},\ \Phi_m(v^m_g)=|G|^{-1/2}\overline{v^m_{\overline g}},\
\Phi_m(v^m_{\overline g})=\tau|G|^{-1/2}\overline{v^m_g},
$$
$$
\Psi_g(\overline{v^g_h})=v^{-g}_{h-g},\ \Psi_g(\overline{v^g_m})=v^{-g}_m,\ \Psi_m(\overline{v^m_g})=|G|^{1/2}v^m_{\overline g},\
\Psi_m(\overline{v^m_{\overline g}})=\tau^{-1}|G|^{1/2}v^m_g,
$$
which implies, for all $g,h,k\in G$:
$$
S(v^g_h\otimes\overline{v^g_k})=v^{-g}_{k-g}\otimes\overline{v^{-g}_{h-g}},\ S(v^g_h\otimes\overline{v^g_m})=v^{-g}_{m}\otimes\overline{v^{-g}_{h-g}},\
S(v^g_m\otimes\overline{v^g_h})=v^{-g}_{h-g}\otimes\overline{v^{-g}_{m}},\
$$
$$
S(v^g_m\otimes\overline{v^g_m})=v^{-g}_{m}\otimes\overline{v^{-g}_m},\
S(v^m_g\otimes\overline{v^m_h})=v^{m}_{\overline h}\otimes\overline{v^{m}_{\overline g}},\
S(v^m_g\otimes\overline{v^m_{\overline h}})=\tau^{-1}(v^{m}_{h}\otimes\overline{v^{m}_{\overline g}}),
$$
$$
S(v^m_{\overline g}\otimes\overline{v^m_h})=\tau(v^{}_{\overline h}\otimes\overline{v^{m}_{g}}),\
S(v^m_{\overline g}\otimes\overline{v^m_{\overline h}})=v^{m}_{h}\otimes\overline{v^{m}_{g}}
$$
and
$$
(v^g_h\otimes\overline{v^g_k})^*=v^{-g}_{h-g}\otimes\overline{v^{-g}_{k-g}},\
(v^g_h\otimes\overline{v^g_m})^*=v^{-g}_{h-g}\otimes\overline{v^{-g}_{m}},\
(v^g_m\otimes\overline{v^g_h})^*=v^{-g}_{m}\otimes\overline{v^{-g}_{h-g}},\
$$
$$
(v^g_m\otimes\overline{v^g_m})^*=v^{-g}_{m}\otimes\overline{v^{-g}_m},\
(v^m_g\otimes\overline{v^m_h})^*=v^{m}_{\overline g}\otimes\overline{v^{m}_{\overline h}},\
(v^m_g\otimes\overline{v^m_{\overline h}})^*=\tau(v^{m}_{\overline g}\otimes\overline{v^{m}_{h}}),
$$
$$
(v^m_{\overline g}\otimes\overline{v^m_h})=\tau^{-1}(v^{m}_{g}\otimes\overline{v^{m}_{\overline h}}),\
(v^m_{\overline g}\otimes\overline{v^m_{\overline h}})=v^{m}_{g}\otimes\overline{v^{m}_{h}}.
$$
We have $B=\underset{g\in G}\oplus B^g\oplus B^m$, where $B^g\cong M_{|G|+1}(\mathbb C),\ \forall g\in G,\ B^m\cong M_{2|G|}(\mathbb C)$


\subsection{Classification of Indecomposable Finite Dimensional $\mathfrak G_{\mathcal T\mathcal Y}$-$C^*$-algebras}

Let us first recall the following well known (see, for instance, \cite{EGNO}, 7.4)

\begin{lemma} \label{point} Equivalence classes of left indecomposable $Vec_G$-module categories are parameterized by couples $(K,\phi)$, where $K$ is a stabilizing subgroup of $G$ and $\phi\in H^2(K,\mathbb C^\times)$.
The set of irreducibles of such a category $\mathcal M(K,\phi)$ is $\Lambda_K=G/K$ and $\phi$ defines the associativities. The corresponding fusion rule is $g\boxtimes\lambda:=g+\lambda,\ \forall g\in G,\ \lambda\in G/K$.

\end{lemma}

Although $\mathcal C=UCorep(\mathfrak G_{\mathcal T\mathcal Y})\cong\mathcal T\mathcal Y(G,\chi,\tau)$, these categories
have different associativities, so we cannot apply directly the classification of module categories from \cite{MeMu}, Section 9, however, we will use similar reasoning. The category $\mathcal C$ is $\mathbb Z_2$-graded, i.e., $\mathcal C=\mathcal C_0\oplus\mathcal C_1$, where $\mathcal C_0\cong Vec_G$ (both these categories have trivial associativities) and $\mathcal C_1$ is generated by a single simple object $U^m$. Indecomposable $\mathcal C_0$-module categories with trivial associativities
are parameterized by their stabilizer subgroups $K<G$, they correspond to $Vec_G$-module categories of the form $\mathcal M(K,{\bf 1})$, where ${\bf 1}$ is the trivial cocycle. Let us denote them by $\mathcal M(K)$.

Then, according to \cite{MeMu}, any indecomposable $\mathcal C$-module category $\mathcal M$ is either indecomposable over $\mathcal C_0$ (we say that it is of type {\bf (I)}, it is then of the form $\mathcal M(K)$) or equivalent to $\mathcal M(K_0)\oplus\mathcal M(K_1)$, where $K_0$ and $K_1$ are subgroups of $G$ (they can be equal) - a category of type {\bf (D)}.

Moreover, $\mathcal C_1$ is an invertible $\mathcal C_0$-bimodule category, so one can define an action of $\mathbb Z_2=<\sigma>$ on the set of (equivalence classes) indecomposable semisimple $\mathcal C_0$-module categories: $\sigma\cdot\mathcal M(K):=\mathcal C_1\boxtimes \mathcal M(K)$.


\begin{notations} \label{perp} For $K<G,\ \rho\in \hat K$ denote $K^\perp_\rho=\{g\in G|\chi(g,k)=\rho(-k),\
\forall k\in K\}$. If $\rho={\bf 1}$ is trivial, denote $K^\perp_{\bf 1}$ by $K^\perp$. Note that $\hat K\cong G/K^\perp$.
\end{notations}

\begin{lemma} \label{orthog}
For any $K<G$, we have $\sigma\cdot\mathcal M(K)=\mathcal M(K^\perp)$.
\end{lemma}
\begin{dm} Adopting the strategy of the proof of \cite{MeMu}, Lemma 30 to our context, let $A_K=\underset{k\in K} \oplus H^k$ be an algebra in the category $\mathcal C_0$ - the analog of the algebra $\mathbb C K$ in $Vec_G$. Viewed as a usual $C^*$-algebra, $A_K$ has the following minimal central orthoprojectors:
$$
P_\lambda=\frac{1}{|K|}\underset{g\in\lambda}\Sigma v^0_{g},\quad P_\rho=\frac{1}{|K|}\underset{k\in K}\Sigma
\rho(k)v^k_m\quad(\lambda\in G/K,\ \rho\in\hat K)
$$
So indecomposable right $A_K$-modules with support in $\mathcal C_0$ are: $V_\lambda=Vec\{v^k_\lambda|k\in K\}$
with the action $v^k_\lambda\circ H^h=v^{h+k}_\lambda\ (h,k\in K)$ and $V_\rho=\mathbb C\underset{k\in K}\Sigma\rho(k)v^k_m$
with the action $(\underset{k\in K}\Sigma\rho(k)v^k_m)\circ H^h=\rho(-h)\underset{k\in K}\Sigma\rho(k)v^k_m\ (h\in H)$, where we denote $v^x_\lambda:=\underset{g\in\lambda}\Sigma v^x_g$ for any $x\in\Omega$. In both cases the stabilizer subgroup is $K$.

Then the category $\mathcal C_1\boxtimes\mathcal M(K)$ can be described as the category of right $A_K$-modules in $\mathcal C$ with support in $\mathcal C_1$ which are of the form $H^m\otimes_R V_\lambda=Vec\{v^m_{\overline p}|p\in\lambda\}$ with the action $v^m_{\overline p}\circ H^h=v^m_{\overline{p+h}}\ (p\in\lambda, h\in K)$ and $H^m\otimes_R V_\rho=Vec\{v^m_r|r\in K^\perp_\rho\}$ with the action $v^m_r\circ H^h=\chi(h,r)v^m_r\ (r\in H^\perp_\rho, h\in H)$.

In order to determine the stabilizer of $H^m\otimes_R V_\lambda$, we calculate, as in the proof of \cite{MeMu}, Lemma 30, for all $g\in G$ the modules $H^g\otimes_R (H^m\otimes_R V_\lambda)=Vec\{\chi(g,p)v^m_{\overline p}|p\in\lambda\}$ with the action $\chi(g,p)v^m_{\overline p}\circ H^h = \chi(g,-h)\chi(g,p+h)v^m_{\overline{p+h}}\ (p\in\lambda, h\in K)$. Therefore, the stabilizer is $K^\perp$.

Similarly, we calculate for all $g\in G$ the modules $H^g\otimes_R (H^m\otimes_R V_\rho)=Vec\{v^m_{r-g}|r\in K^\perp_\rho\}$, but $r-g\in K^\perp_\rho$ is equivalent to $g\in K^\perp$.
\hfill\end{dm}

Thus, in case {\bf (I)} necessarily $K= K^\perp$, so $|G|$ must be a square, and $\Lambda=G/K$. In case {\bf (D)} $\mathcal M\cong\mathcal M(K)\oplus\mathcal M(K^\perp)$ and $\Lambda=G/K\sqcup G/K^\perp$.

\begin{corollary} \label{modGr}

The fusion rules for indecomposable $UCorep(\mathfrak G_{\mathcal T\mathcal Y})$-module categories are: $U^g\boxtimes M_\lambda=M_{g+\lambda}\ (\forall g\in G,\ M_\lambda\in Irr(\mathcal M)$) in all cases and:

$$
\text{For}\ \mathcal M=\mathcal M(K):\ U^m\boxtimes M_\lambda=\underset{\mu\in G/K}\Sigma M_\mu,\ \text{where}\
M_\lambda, M_\mu\in Irr(\mathcal M(K))
$$

For $\mathcal M=\mathcal M(K)\oplus\mathcal M(K^\perp)$:

$$
U^m\boxtimes M_\lambda=\underset{\mu\in G/K^\perp}\Sigma M_\mu,\quad
U^m\boxtimes M_\mu=\underset{\lambda\in G/K}\Sigma M_\lambda,
$$
where $M_\lambda\ (\lambda\in G/K)$ and $M_\mu\ (\mu\in G/K^\perp)$ are in $Irr(\mathcal M)$.
\end{corollary}

\begin{dm} A priori, we have the following fusion rules with $U^m$:

$$
\text{For}\ \mathcal M=\mathcal M(K):\ U^m\boxtimes M_\lambda=\underset{\mu\in G/K}\Sigma n^\mu_\lambda
M_\mu,\ (M_\lambda, M_\mu\in Irr(\mathcal M(K)),\ n^\mu_\lambda\in \mathbb Z_+)
$$

For $\mathcal M=\mathcal M(K)\oplus\mathcal M(K^\perp)$:

$$
U^m\boxtimes M_\lambda=\underset{\mu\in G/K^\perp}\Sigma m^\mu_\lambda M_\mu,\quad
U^m\boxtimes M_\mu=\underset{\lambda\in G/K}\Sigma m^\lambda_\mu M_\lambda,
$$
where $M_\lambda\ (\lambda\in G/K),\ M_\mu\ (\mu\in G/K^\perp)$ are in $Irr(\mathcal M)$ and $m^\mu_\lambda,
m^\lambda_\mu\in \mathbb Z_+$.

The relations of the type $(U^g\otimes U^m)\boxtimes M_\lambda=U^g\boxtimes(U^m\boxtimes M_\lambda)$,
$(U^m\otimes U^g)\boxtimes M_\lambda=U^m\boxtimes(U^g\boxtimes M_\lambda)$ and similar relations with $M_\mu$
show that $n^\mu_\lambda, m^\mu_\lambda$ and $m^\lambda_\mu$ do not depend on $\lambda$ and $\mu$. Then it
remains to apply again $U^m$ to the above equalities and to use the last remark, the relation $U^m\otimes U^m=\underset{g\in G}\Sigma U^g$ and the fact that $|G|=|K||K^\perp|$.\hfill\end{dm}

\begin{corollary} \label{generator}

Any object $M=\underset{\lambda\in\Lambda}\oplus m_{\lambda}M_{\lambda}$ of an indecomposable semisimple
$UCorep(\mathfrak G_{\mathcal T\mathcal Y})$-module category is a generator.
Indeed, Corollary \ref{modGr} shows that already any $\mathcal M_\lambda$ is a generator.

Therefore, the set of all couples $(\mathcal M,M)$ is parameterized:

in case {\bf (I)} by couples $(K,\{m_\lambda|\lambda\in G/K\})$, where $K=K^\perp<G$ and $m_\lambda\in\mathbb Z_{+}$
are such that at least one $m_\lambda > 0$.

in case {\bf (D)} by triples $(K,\{m^0_\lambda|\lambda\in G/K\},\{m^1_\mu|\mu\in G/K^\perp\})$, where $K<G$ and $m^0_\lambda,m^1_{\mu}\in\mathbb Z_{+}$ are such that at least one of them is nonzero.
\end{corollary}

\begin{lemma} \label{isomod} The group $Aut(\mathcal M)$ of autoequivalences of an indecomposable semisimple
$UCorep(\mathfrak G_{\mathcal T\mathcal Y})$-module category $\mathcal M$ with trivial associativities is
as follows:
\vskip 0.5cm
(1) In case {\bf (I)} for any $\phi\in Aut(\mathcal M)$, there exists a unique $p\in G/K$ such that $\phi(M_\lambda)=M_{p+\lambda}$, for all $\lambda \in G/K$, so $Aut(\mathcal M)\cong G/K$.

(2) In case {\bf (D)} and:

\hskip 0.5cm a) $K \neq K^\perp$, for all $\phi \in Aut(\mathcal M)$, there exists a unique $(p_0,p_1)\in G/K\times G/K^\perp$ such that $\phi(M_\lambda)=M_{p_0+\lambda}$ and $\phi(M_\mu)=M_{p_1+\mu}$ for all $\lambda \in G/K, \mu\in
G/K^\perp$, so  $Aut(\mathcal M))\cong G/K\times G/K^\perp$.

\hskip 0.5cm b) $K =K^\perp$,  $Aut(\mathcal M)$, viewed as a bijection of $G/K \times G/K$ on itself, is generated by translations of irreducibles  $(M_\lambda,M_\mu)$ by elements  $(p_0,p_1)\in G/K\times G/K$ and the flip $(M_\lambda,M_\mu) \mapsto (M_\mu,M_\lambda) $. Therefore, $Aut(\mathcal M)\cong (G/K\times G/K)\underset \sigma \ltimes \mathbb Z_2$, where $\sigma$ is the flip of $G/K\times G/K$.


\end{lemma}

\begin{dm} (1) By definition of $\phi$, we must have $\phi(U^g\boxtimes M_\lambda)=U^g\boxtimes\phi(M_\lambda)$,
for all $g\in G,\lambda\in G/K$. Then, putting $M_p=\phi(M_K)$, we have the needed formula for $\phi$. Conversely,
it is easy to check that for such a $\phi$ we have $\phi(U^x\boxtimes M_\lambda)=U^x\boxtimes\phi(M_\lambda)$,
for all $x\in\Omega,\ \lambda\in G/K$.

(2a) As $\mathcal M=\mathcal M(K)\oplus\mathcal M(K^\perp)$ and $\mathcal M(K)\ncong\mathcal M(K^\perp)$, the above result applies to the corresponding restrictions of $\phi$.

(2b) Now the above mentioned components have equal rights, so $\phi$ can permute them and we are done.
\hfill\end{dm}

Corollary \ref{generator} implies that any object $M=\underset{\lambda\in\Lambda}\oplus m_\lambda M_\lambda$ of a module category $\mathcal M$ as above can be identified either with a collection $\{m_\lambda|\lambda\in G/K\}$ or with a double collection $(\{m_\lambda|\lambda\in G/K\},\{m_\mu|\mu\in G/K^\perp\})$, where $m_\lambda,m_\mu\in\mathbb Z_+$.

These considerations and Theorem \ref{main} prove Theorem \ref{TY}.

\begin{remark} \label{dim} Let us compute the dimensions of the spectral subspaces of a finite dimensional
$\mathfrak G_{\mathcal T\mathcal Y}$-$C^*$-algebra $(A,\alpha)$. By Theorem \ref{main}, given a $C^*$-module category $\mathcal M$ over $UCorep(\mathfrak G_{\mathcal T\mathcal Y})$ with a generator $M=\underset{\lambda\in\lambda}\oplus m_\lambda M_\lambda$, we have $A_{U^x}=F(U^x)\otimes \overline H_x\ (\forall x\in\Omega)$, where $F:UCorep(\mathfrak G_{\mathcal T\mathcal Y})\to Corr_f(R)$ is the functor defined by $F(U^x):=Hom(M,U^x\boxtimes M)$, $R=End(M)$. Clearly,
$$
X^x:=F(U^x)\cong\underset{\lambda,\rho\in\lambda}\oplus m_\lambda m_\rho Hom(M_\lambda,U^x\boxtimes M_\rho).
$$
As $Hom(M_\lambda,U^g\boxtimes M_\rho)=\delta_{\lambda,g\cdot\rho}\mathbb C,\ \forall\lambda,\rho\in\lambda$, we have
$dim(X^g)=\underset{\rho\in\lambda}\Sigma m_\rho m_{g\cdot\rho}$.

Now, in case {\bf (I)}, $Hom(M_\lambda,U^m\boxtimes M_\rho)\equiv\mathbb C$, so $dim(X^m)=\underset{\lambda,\rho\in G/K}\Sigma m_\lambda m_{\rho}$.

And in case {\bf (D)}, $Hom(M_\lambda,U^m\boxtimes M_\rho)=0$ when $\lambda,\rho\in G/K$ or $\lambda,\rho\in G/K^\perp$, and
$Hom(M_\lambda,U^x\boxtimes M_\rho)=\mathbb C$ otherwise. So, $dim(X^m)=2\underset{\lambda\in G/K}\Sigma
m_\lambda\times$ $\times\underset{\rho\in G/K^\perp}\Sigma m_\rho$. Therefore, in case {\bf (D)}, $dim X^m$ must be even.
\end{remark}
\end{section}

\begin{section}{Indecomposable Weak Coideals of $\mathfrak G_{TY}$}\label{coid}

We begin the classification of indecomposable weak coideals of $\mathfrak G_{TY}$ by giving a canonical basis for them.

\begin{notations} \label{$v^g_X$} For all $g\in G$ and $X\subset G\sqcup\{m\}$, let us denote:
$$
v^g_X=\underset{x\in X} \Sigma v^g_x,\quad v^m_X=\underset{g\in G\cap X} \Sigma v^m_g,\quad v^m_{\overline X}=
\underset{g\in G\cap X} \Sigma v^m_{\overline g}.
$$
\end{notations}

\begin{lemma} \label{basis} Let $A$ be a weak coideal of $B$. Then:

a) For any $g\in G$ such that $X^g\neq \{0\}\ $, there exists a subset $I^g\subset I^0=\{\Gamma_i|i=1,2,...,k_0\}$ of cardinality
$k_g$ and a set of vectors $\{v^g_{\Gamma_i}(\Theta^g)|i\in I^g\}$ which is a basis of $X^g$, where $\Theta^g$ is a map from $\Gamma^g=\underset{i\in I^g}\sqcup\Gamma_i$ to $\mathbb T:=\{z\in\mathbb C| |z|=1\}$.

b) If $X^m\neq\{0\}$, then $v^0_m\in X^0$, so we can chose $\{m\}\in I^0$, and there exists a subset $I^m\subset I^0\backslash\{m\}$
of cardinality $k_m$ and a basis of $X^m$ of the form $\{v^m_{\Gamma_i}(\Theta^m),v^m_{\overline{\Gamma_i}}(\Theta^m)|i\in I^m\}$,
where $\Theta^m$ is a map from $\Gamma^m=\underset{i\in I^m}\sqcup\Gamma_i$ to $\mathbb T:=\{z\in\mathbb C| |z|=1\}$. If $k_m=k_0-1$,
this weak coideal is indecomposable.
\end{lemma}

{\bf Proof.} a) Let $v^g=\underset{x\in\Omega}\Sigma a_x v^g_x$ be a nonzero vector from $X^g$. Then $v^g=v^g\circ v^0_\Gamma=
\underset{i\in I^0}\Sigma v^g(\Gamma_i)$, where $v^g(\Gamma_i)=\underset{x\in\Gamma_i}\Sigma a_x v^g_x$. Hence $X^g=
\underset{i\in I^g}\oplus X^g_{\Gamma_i}$, where $X^g_{\Gamma_i}\ (i\in I^g)$ are subspaces of $X^g$ containing $v^g(\Gamma_i)
\neq 0$. We have:
$$
v^g(\Gamma_i)^\sharp \circ v^g(\Gamma_i)=\underset{x\in\Gamma_i}\Sigma |a_x|^2 v^0_x=C v^0_{\Gamma_i},\quad\text{where}\quad C>0.
$$

Let $w^g(\Gamma_i)=\underset{y\in\Gamma_i}\Sigma b_y v^g_y\in X^g_{\Gamma_i}$ be another vector with $|b_y|\equiv 1$, then:
$$
\tilde v^g(\Gamma_i)^\sharp \circ w^g(\Gamma_i)=\underset{x\in\Gamma_i}\Sigma \overline{a_x} b_x v^0_x=D v^0_{\Gamma_i},
$$
where $|D|=1$. Then $b_x=D a_x$ for all $x\in\Gamma_i$ which shows that any $X^g_{\Gamma_i} (i\in I^g)$ is generated by a unique,
up to a scalar $D\in\mathbb T$, vector as above. We fix such elements and denote them by $v^g_{\Gamma_i}(\Theta^g)$,
the map $\Theta^g$ being defined by the coefficients of the chosen elements.



b) Let $X^m\neq\{0\}$ and let $v^m=\underset{g\in G}\Sigma a_g v^m_g + \underset{h\in G}\Sigma b_h v^m_{\overline h}$ be its
nonzero vector. Then $(v^m)^\sharp:=\Psi_m(\overline v^m)=|G|^{1/2}(\underset{g\in G}\Sigma\overline a_g v^m_{\overline g} + \underset{h\in G}
\Sigma\overline b_h\tau^{-1} v^m_h$). Next, we compute:
$$
v^m\circ(v^m)^\sharp=|G|^{1/2}(\underset{g,k\in G}\Sigma a_g\overline a_k v^{k-g}_{k} + \underset{p,h\in G}
\Sigma\overline |b_h|^2\overline{\chi(p,h)} v^p_m)
$$
and similarly
$$
(v^m)^\sharp\circ v^m=|G|^{1/2}(\tau\underset{g,p\in G}\Sigma |a_g|^2\overline{\chi(p,g)}v^p_m + \tau^{-1}\underset{h,k\in G}
\Sigma\overline b_k b_h v^{h-k}_{h}).
$$
Hence, the components of index $p$ of these vectors are:
$$
[v^m\circ(v^m)^\sharp]_p=|G|^{1/2}(\underset{g\in G}\Sigma a_g\overline a_{p+g} v^{p}_{p+g} + \underset{h\in G}
\Sigma |b_h|^2\overline{\chi(p,h)} v^p_m)
$$
and similarly
$$
[(v^m)^\sharp\circ v^m]_p=|G|^{1/2}(\tau\underset{g\in G}\Sigma |a_g|^2\overline{\chi(p,g)}v^p_m + \tau^{-1}\underset{k\in G}
\Sigma\overline b_k b_{p+k} v^{p}_{p+k}).
$$
In particular, the components of index $0$ of these vectors are:
$$
[v^m\circ(v^m)^\sharp]_0=|G|^{1/2}(\underset{g\in G}\Sigma |a_g|^2 v^{0}_{g} + (\underset{h\in G}\Sigma |b_h|^2) v^0_m)
$$
and similarly
$$
[(v^m)^\sharp\circ v^m]_0=|G|^{1/2}(\tau(\underset{g\in G}\Sigma |a_g|^2) v^0_m + \tau^{-1}\underset{k\in G}
\Sigma |b_k|^2 v^{0}_{k}).
$$
Since at least one of $a_g$ or $b_h$ is nonzero, it follows that $v^0_m\in X^0$, so we can chose $\{m\}\in I^0$. Further:
$$
v^0_m\circ v^m=\underset{h\in G}\Sigma b_h v^m_{\overline h}\in X^m,\quad v^m\circ v^0_m=\underset{g\in G}\Sigma a_g v^m_g\in X^m,
$$
which shows that $v^0_m\notin Z(A)$ and that $X^m=X^m_1\oplus X^m_2$, where the subspaces $X^m_1, X^m_2\subset X^m$ consist,
respectively, of vectors of the form $\underset{g\in G}\Sigma a_g v^m_g$ and $\underset{h\in G}\Sigma b_h v^m_{\overline h}$. As
$(X^m_1)^\sharp=X^m_2$ and $(X^m_2)^\sharp=X^m_1$, $dim (X^m)$ must be even.

Now, the relations $v^0_{\Gamma_i}\circ v^m=\underset{g\in\Gamma_i}\Sigma a_g v^m_g:=w^m_{\Gamma_i}$ show that $X^m$ has a basis
of the form $\{w^m_{\Gamma_i},(w^m_{\Gamma_i})^\sharp|i\in I^m\}$, and using the same reasoning as in part a), one can normalize:
$w^m_{\Gamma_i}=v^m_{\Gamma_i}(\Theta^m)$. Finally, if $k_m=k_0-1$, there is no a combination of $v^0_{\Gamma_i}$ which would
commute with all $v^m_{\Gamma_i}(\Theta^m)$, so $A$ is indecomposable.\hfill $\square$

Corollary \ref{multiplicity} implies that for weak coideals we have $m_\lambda\in\{0,1\}$ for all $\lambda\in\Lambda$,
so that the generator $M$ can be identified either with a nonempty subset $Z\subset G/K$ or with a couple of subsets
$(Z_0,Z_1)\subset G/K\times G/K^\perp$, at least one of which is nonempty.

\begin{subsection}{The case $A^m=\{0\}$}

\begin{remark} \label{$dim F(U^m)=0$} Let $A$ be an indecomposable weak coideal such that $dim (X^m)=0$. Then either
the set $I^0$ consists of only one subset $\tilde\Gamma\subset\Omega$ containing $\{m\}$ (so that $dim(X^0)=1$) or does
not contain subset $\tilde\Gamma\subset\Omega$ containing $\{m\}$.

Indeed, if $\tilde\Gamma\in I^0$, it suffices to show that $v^0_{\tilde\Gamma}$ commutes with any basis element
$v^g_{\Gamma_i}(\Theta^g)$ of $X$. Using the fact that no more than one element of $I^0$ can contain $\{m\}$ as well as
Lemma \ref{basis}, one can see
that for any $\Gammaà_i\in I^0\backslash\{\tilde\Gamma\}$ we have $v^0_{\tilde\Gamma}\circ v^g_{\Gamma_i}(\Theta^g)=
v^g_{\Gamma_i}(\Theta^g)\circ v^0_{\tilde\Gamma}=0$ and $v^0_{\tilde\Gamma}\circ v^g_{\tilde\Gamma}(\Theta^g)=
v^g_{\tilde\Gamma}(\Theta^g)\circ v^0_{\tilde\Gamma}=v^g_{\tilde\Gamma}(\Theta^g)$. It follows that either the basis of
$X$ consists only from vectors of the form $v^g_{\tilde\Gamma}(\Theta^g)$ or does not contain such vectors at all.
\end{remark}

The equality $dim X^0=1$ implies $M=M_{\lambda_0}$ for some $\lambda_0\in\Lambda$. Then $dim F(U^g)=1$ if
$g\in K$ and $dim F(U^g)=0$ otherwise. This gives a unique, up to isomorphism of $\mathfrak G$-$C^*$algebras,
connected coideal $I^m_K=\underset{k\in K}\oplus(\mathbb C v^k_m\otimes\overline H^k)$.



Now suppose that $\Gamma_i\subset G,\ \forall i\in I^0$. As $dim F(U^m)=0$, $M$ is supported only on $G/K$ or only on $G/K^\perp$. Let us consider the first of these cases, the second one is completely similar. Identify the generator $M$ with a nonempty subset $Z\subset G/K$. The following example shows that any such $Z$ gives rise to an indecomposable weak coideal of $\mathfrak G_{\mathcal T\mathcal Y}$.

\begin{example} \label{gen(iii)}
Let $Z$ be a nonempty subset of $G/K$, then Remark \ref{dim} gives
$dim X^g=\underset{\lambda\in G/K}\Sigma m_\lambda m_{g+\lambda}=|Z\cap(g+Z)|$.

Put $X^g=Vec\{v^g_\lambda|\lambda\in Z\cap(g+Z)\}$ and $X^m=\{0\}$. For any $v^g_\lambda\in X^g\ (g\in G)$, we have $(v^g_\lambda)^\sharp=
v^{-g}_{\lambda-g}\in X^{-g}$. Indeed, as $\lambda\in Z\cap(g+Z)$, there is $\lambda'\in Z$ such that $\lambda=g+\lambda'$,
so $\lambda-g=\lambda'\in Z$. Clearly, $(\lambda-g)\in Z-g$, hence $(\lambda-g)\in Z\cap (Z-g)$. We also have:
$$
v^g_\lambda\circ v^h_\mu=\delta_{\mu,h+\lambda}v^{g+h}_\mu\in X^{g+h}\quad\text{for\ all}\quad v^g_\lambda\in X^g,\ v^g_\mu\in X^h\ (g,h\in G).
$$
Indeed, as $\lambda\in Z\cap(g+Z),\ \mu\in Z\cap(h+Z)$, there are $\lambda',\ \mu'\in Z$ such that $\lambda=g+\lambda',\ \mu=h+\mu'$,
so the above product is nonzero if and only if $\mu=h+\lambda=h+g+\lambda'\in g+h+Z$. Since $\mu\in Z$, it follows that $\mu\in Z\cap (g+h+Z)$.
Thus, Lemma \ref{wcoid}, a) implies that the family $\{X^x|x\in\Omega\}$ generates a weak coideal $A\subset B$ with unit $1_A=v^0_L\otimes\overline v^0_\Omega$, where $L:=\underset{\lambda\in Z}\bigcup\lambda$.
\end{example}

\begin{remark}
$A$ is never a coideal but when $|Z|=1$ it is isomorphic to a connected coideal $I^\Omega_K= \underset{k\in K}\oplus (\mathbb C (\underset{x\in\Omega}\Sigma v^h_x) \otimes\overline H^k)$, which is the "right sided" version
of the left coideal $I_K$ from \cite{M}. $I^\Omega_K$ is also isomorphic to $I^m_K$ above.
If $|Z|>1$, $A$ is also indecomposable because for an arbitrary proper subset $Z_0\subset Z$, the element
$\underset{\lambda\in Z_0}\Sigma v^0_\lambda$ does not commute with any $v^g_\mu\ (\mu\in Z_0,\ g\notin K)$.
\end{remark}

It follows from Remark \ref{dim} that $\mathfrak G_{\mathcal T\mathcal Y}$-$C^*$-algebras $(A,\alpha)$ with
$A^m=\{0\}$ can be only of type (D).

We can summarize the above considerations as follows:

\begin{proposition} \label{clas1} Isomorphism classes of indecomposable weak coideals $A$ of $\mathfrak G_{\mathcal T
\mathcal Y}$ with $A^m=\{0\}$ are parameterized by couples $(K,Z^{orb})$, where $K<G$ and $Z^{orb}$ is the orbit
of a nonempty subset $Z\subset G/K$ or $Z\subset G/K^\perp$ under the action of the group of the translations on $G/K$ (resp., on $G/K^\perp$).
$A$ is isomorphic to a coideal if and only if $|Z|=1$.
\end{proposition}
\end{subsection}


\begin{subsection}{The case $A^m\neq\{0\}$}

\begin{proposition} \label{noexist}
There is no weak coideals of $\mathfrak G_{\mathcal T\mathcal Y}$ corresponding to module categories $\mathcal M$ with $\Lambda=G/K$.
\end{proposition}
\begin{dm}
Let $A$ be such a weak coideal and $M$ be the corresponding generator identified with the subset $Z$ of $G/K$. Then $k_0=dim(X^0)=|Z|$ and $dim(X^m)=|Z|^2$. In terms of Lemma \ref{basis}, b) we have $dim(X^m)=2k_m$, where $k_m\leq k_0-1$, so that $|Z|^2\leq 2(|Z|-1)$ which is only possible if $|Z|=1$. But then $dim(X^m)=1$ - contradicts to the fact that $dim(X^m)$ must be even.
\hfill\end{dm}

\begin{proposition} \label{neq 0}
Let $A$ be a weak coideal of $\mathfrak G_{\mathcal T\mathcal Y}$ corresponding to a module category $\mathcal M$ with $\Lambda=G/K\sqcup G/K^\perp$ and a generator $M$ defined by a nonempty subset $(Z_0,Z_1)\subset G/K\times G/K^\perp$.
Then either $|Z_0|=1$ or $|Z_1|=1$.
\end{proposition}
\begin{dm}
We have $k_0=dim(X^0)=|Z_0|+|Z_1|$ and $dim(X^m)=2|Z_0||Z_1|$. In terms of Lemma \ref{basis}, b) we have $dim(X^m)=2k_m$, where $k_m\leq k_0-1$, so $|Z_0||Z_1|\leq |Z_0|+|Z_1|-1$ from where either $|Z_0|=1$ or $|Z_1|=1$.
\hfill\end{dm}

The following example shows that any such set $(Z_0,Z_1)$ gives rise to an indecomposable weak coideal of $\mathfrak G_{\mathcal T\mathcal Y}$.

\begin{example} \label{new'}
Let $Z$ be a nonempty subset of $G/K$ and $\rho_0\in G/K^\perp$. For the generator corresponding to $Z\sqcup \rho_0$ we have $dim X^m=2|Z|$, $dim X^g=|Z\cap(g+Z)|$ if $g\notin K^\perp$ and $dim X^g=|Z\cap(g+Z)|+1$ if $g\in K^\perp$.

Put $X^m=Vec\{v^m_\lambda, v^m_{\overline \mu}|\lambda,\mu\in Z\}$, $X^g=Vec\{v^g_m, v^g_\lambda|\lambda\in Z\cap(g+Z)\}$ if $g\in K^\perp$
and $X^g=Vec\{v^g_\lambda|\lambda\in Z\cap(g+Z=\}$ if $g\notin K^\perp$. The next relations, where $g,h\in G,k,l\in K^\perp, \lambda,\mu\in Z, u(\lambda)$ is a representative of the coset $\lambda$, show that the family $\{X^x|x\in\Omega\}$ satisfies the conditions a) of Lemma \ref{wcoid}:
$$
(v^m_\lambda)^\sharp=|G|^{1/2}v^m_{\overline \lambda},\ (v^m_{\overline \lambda})^\sharp=\tau^{-1}|G|^{1/2}v^m_\lambda,\
(v^k_m)^\sharp=v^{-k}_m,\ (v^g_\lambda)^\sharp=v^{-g}_{\lambda-g},
$$
$$
v^k_m\circ v^l_m=v^{k+l}_m,\ v^g_\lambda\circ v^h_\mu=\delta_{\mu,g+\lambda}v^{g+h}_\mu,\ v^k_m\circ v^h_\lambda=v^h_\lambda\circ v^k_m=v^k_m\circ v^m_\lambda=0,
$$
$$\ v^k_m\circ v^m_{\overline\lambda}=\chi(k,u(\lambda))v^m_{\overline\lambda}, v^m_{\overline\lambda}\circ v^k_m=0,\ v^m_\lambda\circ v^k_m=\chi(u(\lambda),k)v^m_\lambda,
$$
$$
v^g_\lambda\circ v^m_\mu=\delta_{\lambda,\mu}v^m_{\lambda-g},\ v^m_{\overline\lambda}\circ v^g_\mu=\delta_{\lambda,\mu}v^m_{\overline{g+\lambda}},
\ v^m_\lambda\circ v^g_\mu=v^g_\lambda\circ v^m_{\overline\mu}=v^m_\lambda\circ v^m_\mu=v^m_{\overline\lambda}\circ v^m_{\overline\mu}=0,
$$
and finally, using the fact that $\underset{k\in K}\Sigma\chi(g,k)=|K|$ if $g\in K^\perp$ and is $0$ otherwise:
$$
v^m_\lambda\circ v^m_{\overline\mu}=\underset{g\in(\mu-\lambda)}\Sigma v^g_\mu,\ v^m_{\overline\lambda}\circ v^m_\mu=
\tau |K|\delta_{\lambda,\mu}\underset{k\in K^\perp}\Sigma v^k_m.
$$
So, this family generates an indecomposable weak coideal $A\subset B$, $1_A=(v^0_m + v^0_L)\otimes\overline v^0_\Omega$,
where $L=\underset{\lambda\in Z}\bigcup\lambda$. $A$ is a coideal if and only if $L=G$ in which case it is the analog of
the left connected coideal $J_K$ constructed in \cite{M}.
\end{example}

Now we can summarize the above considerations as follows:

\begin{proposition} \label{clas2} Isomorphism classes of indecomposable weak coideals $A$ of $\mathfrak G_{\mathcal T\mathcal Y}$ with $A^m\neq\{0\}$ are parameterized by pairs $(K,(Z_0,Z_1)^{orb})$, where $K<G$ and $(Z_0,Z_1)^{orb}$ is the orbit of a subset $(Z_0,Z_1)\subset G/K\times G/K^\perp$ such that $min\{|Z_0|,|Z_1|\}=1$ under the action of:

a) the group $G/K\times G/K^\perp$ by  translations, if $K\neq K^\perp$;

b) the  semi direct  product   $(G/K\times G/K)\underset \sigma \ltimes\mathbb Z_2$ generated by  the group $G/K\times G/K$ acting by translations and the  flip  $\sigma: (Z_0,Z_1) \mapsto  (Z_1,Z_0)$ if $K=K^\perp$.

$A$ is isomorphic to a coideal if and only if either $Z_0=G/K$ or $Z_1=G/K^\perp$.
\end{proposition}

Finally, Theorem \ref{clas3} follows from Propositions \ref{clas1} and \ref{clas2}.




\end{subsection}
\end{section}
\bibliographystyle{plain}
\bibliography{biblio1}

\end{document}